\documentclass{article}
\usepackage{graphicx}
\usepackage{amsmath}
\usepackage{booktabs}
\usepackage{caption}
\usepackage{geometry}
\graphicspath{{figures/}}
\begin{document}

\title{ A \ time-space B-spline integrator for the Burgers' equation }
\author{ Idris Dag$^{1}$, Serkan U\u{g}urluo\u{g}lu$^{2}$ and Nihat Adar$^{1}$ \\
Faculty of Engineering and Architecture, Computer Engineering Department$%
^{1} $\\
Faculty of Education $^{2}$, \\
Eskisehir Osmangazi University$,$T\"{u}rkiye$.$ }
\date{}
\maketitle

\begin{abstract}
The purpose of this paper is to propose a new algorithm for obtaining
approximate solutions to the Burgers' equation (BE). Integration in time by
a quadratic B-spline collocation method is shown. To the best of our
knowledge, B-splines have not previously been used to integrate partial
differential equations in both time and space. First, the BE is integrated
using quadratic B-spline functions in time, and then the time-integrated BE
is further solved in space via the cubic B-spline collocation method. The
resulting recursive algebraic equation is used to obtain both shock wave and
front propagation solutions of the BE, demonstrating the effectiveness of
the space--time B-spline collocation method.

keywords: Burgers' equation, collocation method, spline, numerical solutions
of partial differential equations
\end{abstract}

\section{Introduction}

Partial differential equations (PDEs) have been solved numerically in
parallel with advancements in computing, and their numerical treatment
typically leads to systems of algebraic equations. This is accomplished by
applying time and space integrators in sequence. Commonly used time
integrators include the Crank--Nicolson scheme, Runge--Kutta methods\cite{mt}%
, geometric integrators\cite{gi1}, exponential integrators\cite{ei,ei1}, and
explicit--implicit methods\cite{nb}. Studies have shown that higher-order
time integrators provide highly accurate solutions to PDEs.

Some of these methods yield discrete solutions, while others provide
functional solutions over the domain of the PDE. A large number of papers
have been published on obtaining numerical solutions to PDEs using these
approaches. Variants of the collocation \cite{cm,aa,st}, Galerkin\cite{as,bi}%
, meshless \cite{ms}, and differential quadrature\cite{dq} methods have been
proposed to solve PDEs numerically.

Nonlinear PDEs generally do not admit analytical solutions, except in
special cases. To explore the solution space of PDEs more thoroughly,
numerical methods are employed using efficient computational techniques. One
such PDE is the Burgers' equation, which models various phenomena in
engineering sciences such as turbulence, gas dynamics, heat conduction,
elasticity, traffic flow, and acoustic waves. The Burgers' equation provides
a mathematical formulation that describes the interplay between convection
and diffusion. Although a series solution exists for the Burgers' equation,
it becomes impractical for small values of the viscosity constant due to
slow convergence. Therefore, numerical methods have been developed to obtain
sharper solutions of the Burgers' equation for smaller viscosity values. The
numerical methods are used with B-splines, which are piecewise polynomials
that are easy to handle in mathematical operations such as integration,
differentiation, etc. The terms of the PDE, that involve higher-order
derivatives, can be approximated by choosing B-splines of the appropriate
degree to satisfy the continuity conditions. Thus, many methods based on
splines have been employed to compute accurate numerical solutions to the
Burgers' equation.

Bateman \cite{ba}, first, suggested the quasi linear parabolic differential
equation known as the one dimensional time dependent Burgers' equation. 
\begin{equation}
{u_{t}}+u{u_{x}}-\upsilon {u}_{xx}=0,\;x\in \lbrack a,b],\;t\in (0,T]
\label{1}
\end{equation}%
along with the initial condition (IC) and boundary conditions (BCs): 
\begin{equation}
u(x,0)=f(x),  \label{2}
\end{equation}%
\begin{equation}
u(a,t)=u_{0},\text{ }u_{x}(a,t)=u_{0}^{\prime },\quad u(b,t)=u_{m},\text{ }%
u_{x}(b,t)=u_{m}^{\prime },  \label{3}
\end{equation}%
where $\upsilon $ is viscosity constant. Small values of $\upsilon $ cause
the formation of shock waves, so both analytical and numerical methods may
not provide reasonable solutions. The unknown $u$ represents the quantity
evolving over time. The full discretization of Burgers' equation is carried
out using quadratic B-splines for the time domain and cubic B-splines for
the spatial domain. To the best of our knowledge, this is the first study in
which B-splines have been employed simultaneously as time and space
integrators for the Burger's equation. The robustness of the proposed
time--space B-spline collocation method is demonstrated through simulations
of both shock and front wave propagation.

\section{Numerical approach}

Consider equally dividing the time domain $[c,d]$ at the time-grid points
such that a mesh $c=t^{0}<t^{1}<...<t^{k}=d$, $t^{j}=t^{0}+j\Delta t$, $%
j=1,....n$ and together with\ artificial time-grid points $%
t^{-2},t^{-1},t^{n+1},t^{n+2}$ outside the time domain $[c,d].$

The quadratic B-spline functions, $B^{j}(t)\in C^{1}[a,b],$ $j=-1,...,n,$ at
the time-grid points in the interval $[c,d]$ are defined as

\begin{equation}
\begin{tabular}{l}
$B^{j}(t)=\dfrac{1}{\Delta t}\left \{ 
\begin{tabular}{ll}
$(t^{j+2}-t)^{2}-3(t^{j+1}-t)^{2}+3(t^{j}-t)^{2},$ & $t^{j-1}\leq t\leq
t^{j} $ \\ 
$(t^{j+2}-t)^{2}-3(t^{j+1}-t)^{2},$ & $t^{j}\leq t\leq t^{j+1}$ \\ 
$(t^{i+2}-t)^{2},$ & $t^{j+1}\leq t\leq t^{j+2}$ \\ 
0 & otherwise%
\end{tabular}%
\right. $%
\end{tabular}
\label{4}
\end{equation}

The analytical solution $u(x,t)$ of the BE is approximated as an expansion
of the B-spline basis functions

\begin{equation}
U(x,t)=\sum_{j=-1}^{n}\delta (x,t^{j})B^{j}(t)  \label{5}
\end{equation}%
where the variables $\delta (x,t^{j})=\delta ^{j}.$

Values of $B^{j}(t)$ and $(B^{j}(t))^{\prime }$ at the grid points are
worked out using B-splines and its first derivatives \ref{4} as

\[
\begin{tabular}{l}
Table 1 \\ 
$%
\begin{tabular}{lllll}
& $t^{j-2}$ & $t^{j-1}$ & $t^{j}$ & $t^{j+1}$ \\ 
$B^{j}(t)$ & $0$ & $1$ & $1$ & $0$ \\ 
$(B^{j}(t))^{\prime }$ & $0$ & $-2/\Delta t$ & $2/\Delta t$ & $0$%
\end{tabular}%
$%
\end{tabular}%
\]

The approximate solution \ref{5} and its derivative at the grid points can
be given in the form of the coefficients $\delta ^{j}$:

\begin{equation}
\begin{tabular}{l}
$U(x,t^{j})=U^{j}=\delta ^{j-1}+\delta ^{j},j=0,...n$ \\ 
$U_{t}(x,t^{j})=(U_{t})^{j}=\dfrac{2}{\Delta t}(-\delta ^{j-1}+\delta ^{j})$%
\end{tabular}
\label{6}
\end{equation}

\bigskip\ \qquad Consider equally dividing the time domain $[a,b]$ at the
time-grid points such that a mesh $a=x_{0}<x_{1}<...<x_{m}=b$, $%
x_{i}=x_{0}+ih$, $x=1,....m$ and together with \ artificial time-grid points 
$x_{-3},x_{-2},x_{-1},$ $x_{m+1},x_{m+2},x_{m+3}$ outside the time domain $%
[a,b].$

We consider \ the cubic B-spline functions $B_{i}(x)\in C^{2}[a,b],$ $%
i=-1,...,m$ at the space -grid points in the interval $[a,b]$: 
\begin{equation}
B_{i}\left( x\right) =\frac{1}{h^{3}}\left\{ 
\begin{array}{ll}
(x-x_{i-2})^{3}, & \text{if }x\in \left[ x_{i-2},x_{i-1}\right] \\ 
{h^{3}}+3{h^{2}}(x-{x_{i-1}})+3h{(x-{x{_{i-1}}})^{2}}-3{(x-{x_{i-1}})^{3}},
& \text{if }x\in \left[ x_{i-1},x_{i}\right] \\ 
{h^{3}}+3{h^{2}}({x_{i-1}}-x)+3h{({x_{i-1}}-x)^{2}}-3{({x_{i-1}}-x)^{3}}, & 
\text{if }x\in \left[ x_{i},x_{i+1}\right] \\ 
(x_{i+2}-x)^{3}, & \text{if }x\in \left[ x_{i+1},x_{i+2}\right] \\ 
0, & \text{otherwise}.%
\end{array}%
\right.  \label{7}
\end{equation}

\begin{equation}
\delta (x,t^{j})=\sum_{i=-1}^{m+1}\sigma _{i}^{j}B_{i}(x)  \label{8}
\end{equation}%
$B_{i}(x)$ and $B_{i}^{^{\prime }}(x)$ are evaluated at the space grid
points seen in table 2:

\[
\begin{tabular}{l}
Table 2 \\ 
$%
\begin{tabular}{llllll}
& $x_{i-2}$ & $x_{i-1}$ & $x_{i}$ & $x_{i+1}$ & $x_{i+2}$ \\ 
$B_{i}(t)$ & $0$ & $1$ & $4$ & $1$ & $0$ \\ 
$B_{i}^{\prime }(t)$ & $0$ & $-3/h$ & $0$ & $3/h$ & $0$ \\ 
$B_{i}^{\prime \prime }(t)$ & $0$ & $6/h^{2}$ & $-12/h^{2}$ & $6/h^{2}$ & $0$%
\end{tabular}%
$%
\end{tabular}%
\]

An approximation and its derivative at the space-grid points can be
expressed in terms of the coefficients $\delta _{i}^{j}$:

\begin{equation}
\begin{tabular}{l}
$\delta (x_{i},t^{j})=\delta _{i}^{j}=\sigma _{i-1}^{j}+4\sigma
_{i}^{j}+\sigma _{i+1}^{j},j=0,...n,i=0,...m$ \\ 
$\delta _{x}(x_{i},t^{j})=(\delta _{x})_{i}^{j}=\dfrac{3}{h}(-\sigma
_{i-1}^{j-1}+\sigma _{i+1}^{j}),$ \\ 
$\delta _{xx}(x_{i},t^{j})=(\delta _{xx})_{i}^{j}=\dfrac{6}{h^{2}}(\sigma
_{i-1}^{j}-2\sigma _{i}^{j}+\sigma _{i+1}^{j}),$%
\end{tabular}
\label{8a}
\end{equation}%
To start the time recursion, time discretization is achieved by substituting
equation \ref{8a} into the BE \ref{1}:

\begin{equation}
\begin{tabular}{l}
$\dfrac{2}{\Delta t}(-\delta ^{j-1}+\delta ^{j})+(\delta ^{j-1}+\delta
^{j})(\delta _{x}^{j-1}+\delta _{x}^{j})-\upsilon (\delta _{xx}^{j-1}+\delta
_{xx}^{j})=0,$%
\end{tabular}
\label{8b}
\end{equation}

The approximate solution $\delta (x,t^{j})$ at the space-grid points is
given by an expansion of the cubic B-splines \ref{8a} as: 
\begin{equation}
\delta ^{j}\left( {x}_{i}{,t}^{j}\right) =\sum_{i=-1}^{m+1}\sigma {_{i}^{j}}%
\left( t\right) {B_{i}}\left( x\right) ,  \label{9}
\end{equation}%
$\sigma {_{i}^{j}}\left( t\right) $ are time dependent parameters determined
by the collocation method. The approximate solution and its first two
derivatives at the grid points can be computed from the cubic B-spline
functions as: 
\begin{equation}
\begin{tabular}{l}
$\delta _{i}^{j}=\delta \left( {x}_{i}{,t}^{j}\right) =\sigma
_{i-1}^{j}+4\sigma _{i}^{j}+\sigma _{i+1}^{j},$ \\ 
$(\delta ^{^{\prime }})_{i}^{j}=\delta ^{\prime }(x_{i},t^{j})=\dfrac{3}{h}%
(-\sigma _{i-1}^{j}+\sigma _{i+1}^{j})$, \\ 
$(\delta ^{\prime \prime })_{i}^{j}=\delta ^{\prime \prime }(x_{i},t^{j})=%
\dfrac{6}{h^{2}}(\sigma _{i-1}^{j}-2\sigma _{i}^{j}+\sigma _{i+1}^{j})$.%
\end{tabular}
\label{11}
\end{equation}%
Replacing approximations of $\delta _{i}^{j},(\delta _{i}^{^{j}})^{\prime }$
and $(\delta _{i}^{j})^{\prime \prime }$ \ref{11} at the space-grid points
in time-discretized equation~\ref{8b} yields $n+1$ the algebraic equations
with $n+3$ unknown parameters:

\begin{equation}
\begin{tabular}{l}
$(1+\dfrac{\Delta t}{2}\delta _{x}^{j-1}+\dfrac{\Delta t}{2}\delta _{x}^{j}-%
\dfrac{3\Delta t}{2h}\delta ^{j-1}-\dfrac{3\upsilon \Delta t}{2h^{2}})\sigma
_{i-1}^{j}+(4+2\Delta t\delta _{x}^{j-1}+2\Delta t\delta _{x}^{j}+\dfrac{%
6\upsilon \Delta t}{h^{2}})\sigma _{i}^{j}$ \\ 
$+(1+\dfrac{\Delta t}{2}\delta _{x}^{j-1}+\dfrac{\Delta t}{2}\delta _{x}^{j}+%
\dfrac{3\Delta t}{2h}\delta ^{j-1}-\dfrac{3\upsilon \Delta t}{2h^{2}})\sigma
_{i+1}^{j}=\delta ^{j-1}-\dfrac{\Delta t}{2}\delta ^{j-1}\delta _{x}^{j-1}+%
\dfrac{\upsilon \Delta t}{2}\delta _{xx}^{j-1}$%
\end{tabular}
\label{12}
\end{equation}%
where $i=0,...,m$

The boundary conditions ,$u(a,t^{j})=U_{0},$and $u(b,t^{j})=U_{m}$ at both
ends, are used to have equations 
\begin{equation}
\sigma _{0}^{j}+4\sigma _{1}^{j}+\sigma _{2}^{j}=U_{0},\sigma
_{m-1}^{j}+4\sigma _{m}^{j}+\sigma _{m+1}^{j}=U_{m}  \label{13}
\end{equation}

Before solving the nonlinear system \ref{12} together with equations\ref{13}%
, the system must be linearized by choosing values from the previous time
step for the terms in equations \ref{12} $\dfrac{\Delta t}{2}\delta
_{x}^{j}\rightarrow \dfrac{\Delta t}{2}\delta _{x}^{j-1}$ and then applying
the following iteration process two or three times at each time step

\begin{equation}
\sigma _{i}^{\ast (j+1)}=\sigma _{i}^{j}+\dfrac{1}{2}(\sigma
_{i}^{j+1}-\sigma _{i}^{j}),\text{ }i=-1,...,m+1  \label{14}
\end{equation}

Initial parameter vector $\mathbf{\sigma }^{0}=(\sigma _{-1}^{0}\mathbf{,}%
\sigma _{0}^{0}\mathbf{,...}\sigma _{m+1}^{0})$ is needed to start the
system given equations \ref{13}. Initial condition and space derivatives at
both ends lead to a system of equation for the vector $\mathbf{\sigma }^{0}:$

\[
\begin{tabular}{l}
$\dfrac{3}{h}(-\sigma _{-1}^{0}+\sigma _{1}^{0})=U_{x}(x_{0},t^{0}),$ \\ 
$\sigma _{i-1}^{0}+4\sigma _{i}^{0}+\sigma
_{i+1}^{0}=U(x_{i},t^{0}),i=0,...,m$ \\ 
$\dfrac{3}{h}(-\sigma _{m-1}^{0}+\sigma _{m+1}^{0})=U_{x}(x_{m},t^{0}).$%
\end{tabular}%
\]

\subsection{Numerical experiments}

The discrete $L_{2}$ and $L_{\infty }$ error norms are 
\[
\begin{tabular}{ll}
$\left\vert U-U_{N}\right\vert _{\infty }=$max$\left\vert
U_{i}^{j}-(U_{N})_{i}^{j}\right\vert $ & $\left\vert U-U_{N}\right\vert _{2}=%
\sqrt{\sum_{i=}^{n}\left\vert U_{i}^{j}-(U_{N})_{i}^{j}\right\vert }$%
\end{tabular}%
\]%
computed at some times to show the efficiency of the method.

a)The burgers' equation provides a shock-like solution expressed as

\begin{equation}
u(x,t)=\dfrac{x/t}{1+\sqrt{t/t_{0}}\exp (x^{2}/(4\upsilon t))},t\geq 1,0\leq
x\leq 1  \label{15}
\end{equation}%
where $t_{0}=exp(1/8\upsilon ).$

The initial condition $U(x,1)$ is used to calculate initial parameters, and
the boundary conditions $u(0,t)=u(1,t)=0$ are incorporated to the numerical
algorithm to control numerical solution at both ends. The algorithm is run
up to time $t=3.25$ with the data set $h=\Delta t=0.01,0.05$ and $0.001$ for 
$\upsilon =0.01$ and $0.005$ over the interval $[0,1]$. The maximum norm, $%
L_{2}$- norm, and peak value of the shock are tabulated in Table 4. Shock
behavior is visualized in Figs. 1 and 3. From these figures, the initial
shock can be seen to shrink over time. The steeper the initial shock used,
the sharper the resulting shock. The magnitude of the steepness decreases as
time progresses. The maximum norm are depicted in Figs. 2 and 4, where
highest errors are observed around the peak of the shock.

	\begin{figure}[h!]
	\centering
	\begin{tabular}{ll}
		\includegraphics[width=0.45\textwidth]{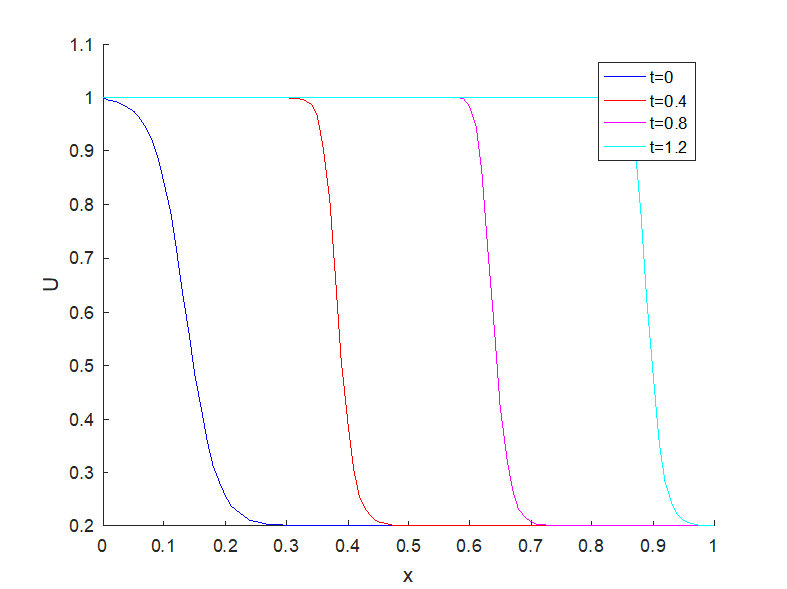} &
		\includegraphics[width=0.45\textwidth]{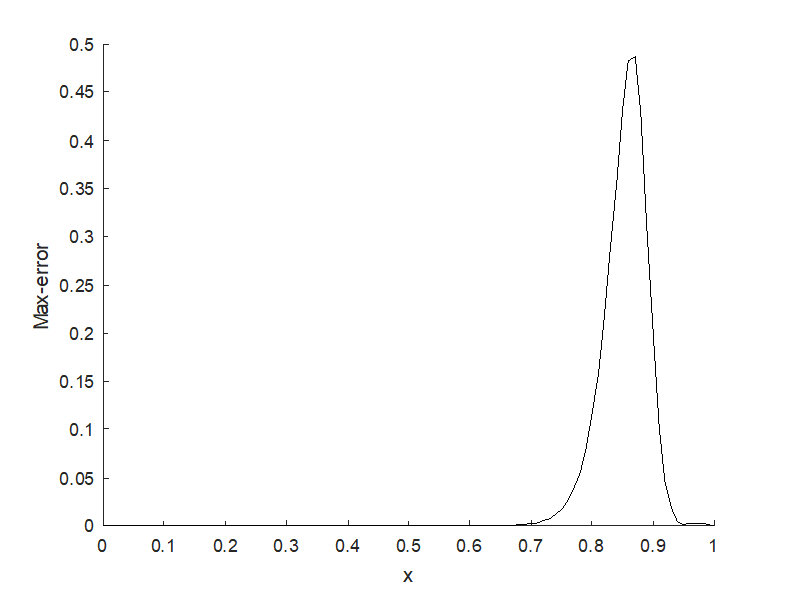} \\
		\textbf{Fig. 6:} Solutions. &
		\textbf{Fig. 7:} Error at time $t = 1.2$.\\
		
		\includegraphics[width=0.45\textwidth]{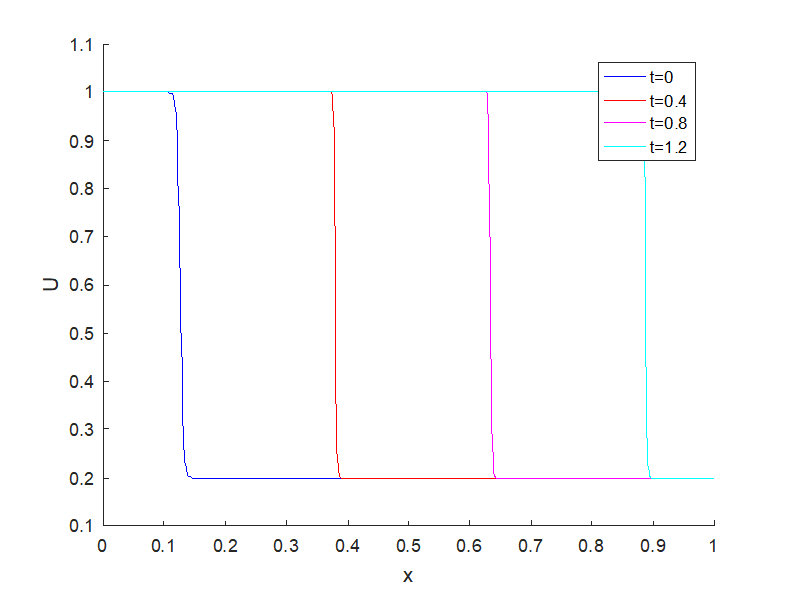} &
		\includegraphics[width=0.45\textwidth]{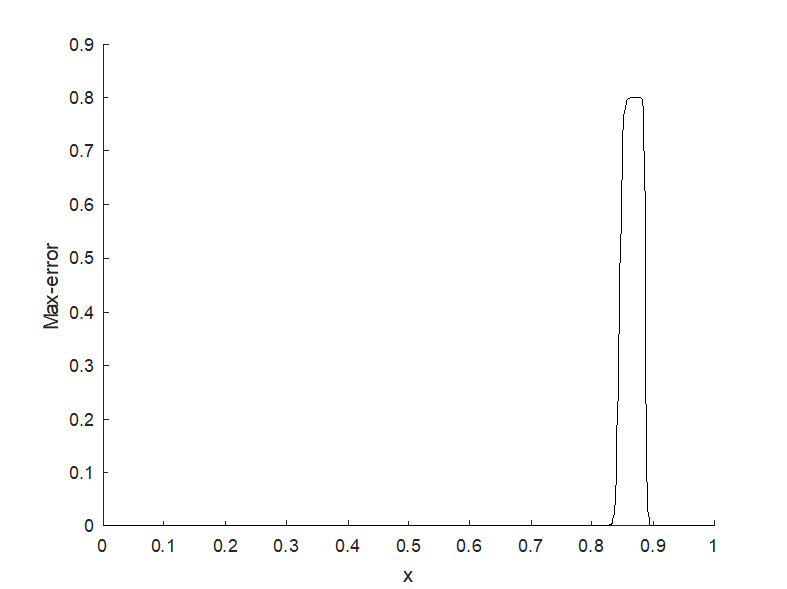} \\
		\textbf{Fig. 8:} Solutions. &
		\textbf{Fig. 9:} Error at time $t = 1.2$.
	\end{tabular}
\end{figure}

\[
\begin{tabular}{l}
Table 4 : Errors and shock peak values. \\ 
\begin{tabular}{ll}
\begin{tabular}{l}
$\ \upsilon $=0.01 \\ 
\begin{tabular}{cccc}
$h$=$\Delta t$ & $L_{\infty }$ & $L_{2}$ & Peak value \\ 
\multicolumn{1}{l}{0.01} & \multicolumn{1}{l}{0.0498} & \multicolumn{1}{l}{
0.0211} & \multicolumn{1}{l}{0.2233} \\ 
\multicolumn{1}{l}{0.005} & \multicolumn{1}{l}{0.0508} & \multicolumn{1}{l}{
0.0207} & \multicolumn{1}{l}{0.2225} \\ 
\multicolumn{1}{l}{0.001} & \multicolumn{1}{l}{0.0516} & \multicolumn{1}{l}{
0.0206} & \multicolumn{1}{l}{0.2203}%
\end{tabular}%
\end{tabular}
& 
\begin{tabular}{l}
$\ \upsilon $=0.005 \\ 
\begin{tabular}{ccc}
$L_{\infty }$ & $L_{2}$ & Peak value \\ 
\multicolumn{1}{l}{0.0582} & \multicolumn{1}{l}{0.0166} & \multicolumn{1}{l}{
0.2503} \\ 
\multicolumn{1}{l}{0.0571} & \multicolumn{1}{l}{0.0159} & \multicolumn{1}{l}{
0.2480} \\ 
\multicolumn{1}{l}{0.0560} & \multicolumn{1}{l}{0.0155} & \multicolumn{1}{l}{
0.2467}%
\end{tabular}%
\end{tabular}%
\end{tabular}%
\end{tabular}%
\]

b) The travelling front wave simulation of the Burgers' equation has
solution 
\begin{equation}
u(x,t)=\dfrac{\alpha +\mu +(\mu -\alpha )\exp (\eta )}{1+\exp (\eta )}
\label{16}
\end{equation}%
where $\eta =\dfrac{\alpha (x-\mu t-\gamma )}{\upsilon },\alpha ,\mu $ and $%
\gamma $ are constants. Boundary conditions \ are $u(0,t)=1$ and $%
u(1,t)=0.2,t\geq 0.$ The initial condition is extracted from the analytical
solution \ref{13} at $t=0$,initially situated at $x=$ $\gamma ,$ propagating
to the right with speed $\mu .$ Parameters $\alpha =0.4,\mu =0.6$ , $\gamma
=0.125,\ $ used \ in the studies are chosen from previous works \cite{cm, bi}%
. Calculations are performed with the space-time increments $\Delta t=$ $%
h=0.01,0.005\ $and $0.001$ for $\upsilon =0.01$ and $0.001$ recursively$.$
Visual solutions \ are depicted using \ the time-space increment $\Delta
t=h=0.01$ for $0.01$ and $0.005$ in figs 6 and 7. Figures 7 and 9 display
the maximum errors at $t=1.2$

The smaller viscosity constant causes sharper wave behavior and higher
error. Larger errors are observed near front of wave \ across the right
boundary. $L_{\infty }$-maximum norm and $L_{2}$-norm \ are documented in
table 5. \ Smaller time-space increments are used,slightly better accuracy
are obtained. But errors for this test problem are larger than the first
test problem.

	\begin{figure}[h!]
	\centering
	\begin{tabular}{ll}
		\includegraphics[width=0.45\textwidth]{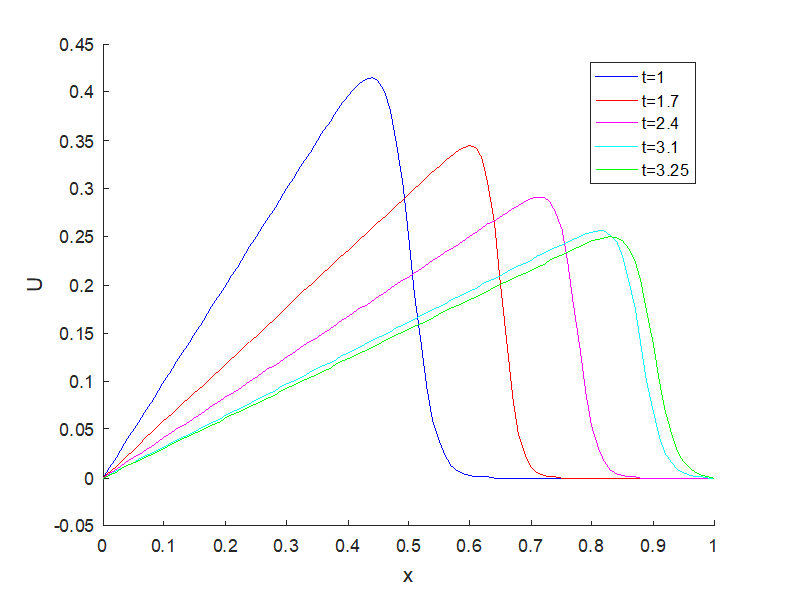} &
		\includegraphics[width=0.45\textwidth]{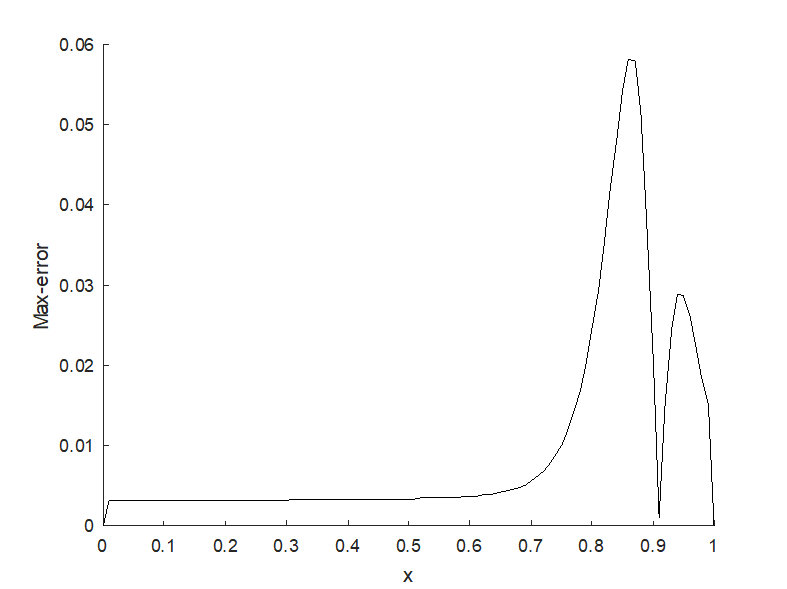} \\
		\textbf{Fig. 1:} Solutions. &
		\textbf{Fig. 2:} Error at time $t = 1.2$.\\
		
		\includegraphics[width=0.45\textwidth]{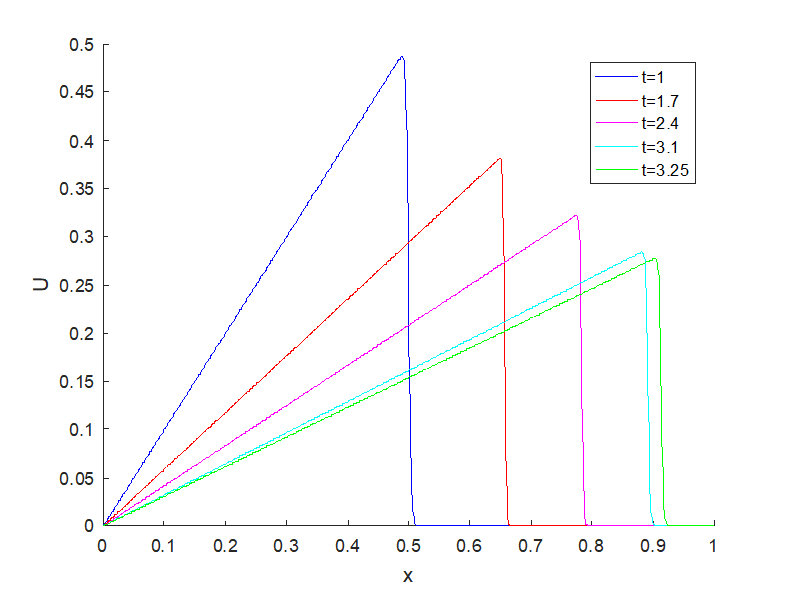} &
		\includegraphics[width=0.45\textwidth]{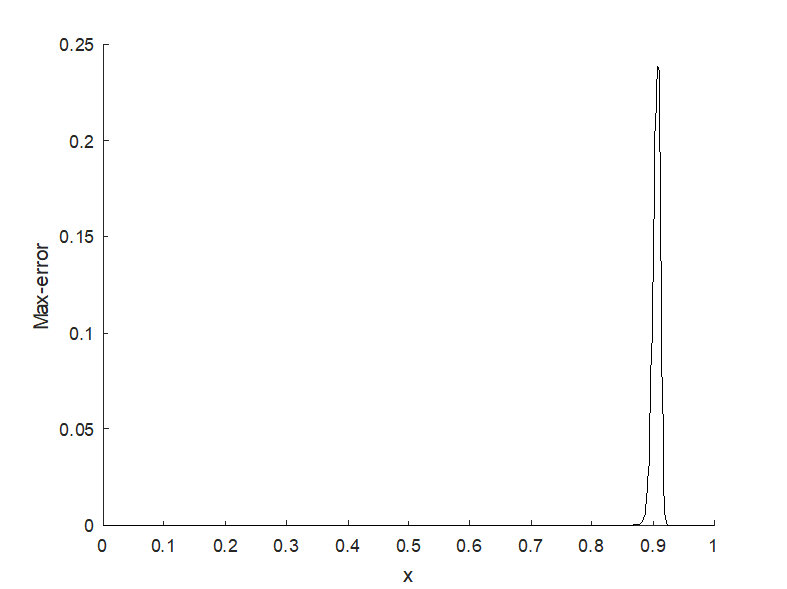} \\
		\textbf{Fig. 3:} Solutions. &
		\textbf{Fig. 4:} Error at time $t = 1.2$.
	\end{tabular}
\end{figure}


	\[
	\begin{tabular}{l}
		Table 5 : Errors \\ 
		\begin{tabular}{ll}
			\begin{tabular}{l}
				$\upsilon $=0.01 \\ 
				\begin{tabular}{lll}
					$h$=$\Delta t$ & $L_{\infty }$ & $L_{2}$ \\ 
					0.01 & 0.4875 & 0.1138 \\ 
					0.005 & 0.2964 & 0.0638 \\ 
					0.001 & 0.1514 & 0.0354%
				\end{tabular}%
			\end{tabular}
			& 
			\begin{tabular}{l}
				$\upsilon $=0.005 \\ 
				\begin{tabular}{ll}
					$L_{\infty }$ & $L_{2}$ \\ 
					0.8841 & 0.2532 \\ 
					0.6906 & 0.1329 \\ 
					0.2427 & 0.0367%
				\end{tabular}%
			\end{tabular}%
		\end{tabular}%
	\end{tabular}%
	\]

\section{Conclusion}

An simple method is introduced for solving the BE and second aspect of the
paper is also to propose an alternative numerical method for integrating the
BE. Both shock and front wave propagations are simulated accurately over
time. Using smaller time--space increments yields smoother solutions. Future
work will focus on improving accuracy by employing higher-order B-splines
for the time--space integration of the Burgers' equation.

\end{document}